\documentclass[12pt]{article}
\usepackage{mathrsfs}
\usepackage{epic,eepic,epsf,epsfig}
\usepackage{amsfonts,srcltx,mathrsfs}
\usepackage{multirow}
\usepackage{amsmath}
\usepackage{amssymb}
\usepackage{amsbsy}
\usepackage{graphicx}
\usepackage{amsfonts}
\usepackage{color}
\usepackage{setspace}

\newtheorem{lem}{Lemma}[section]%
\newtheorem{theorem}[lem]{Theorem}%
\newtheorem{prop}[lem]{Proposition}%

\def\a{\alpha} \def\b{\beta}  \def\d{\delta}

 \def\O{\Omega} \def\G{\Gamma}

\def\di{\bigm|} \def\lg{\langle} \def\rg{\rangle}
\def\nd{\mathrel{\bigm|\kern-.7em/}}

\def\f{\noindent}

\def\P\GammaL{\hbox{\rm P\GammaL}}

\def\Aut{\hbox{\rm Aut}}

\def\mod{\hbox{\rm mod }}

\newcommand{\qed}{\mbox{\raisebox{0.7ex}{\fbox{}}} \vspace{4truemm}}
\def\mz{{\mathbb Z}}

\begin{document}
\title{On cubic symmetric non-Cayley graphs with solvable
automorphism groups}

\footnotetext[1]{ Corresponding author. E-mails:
yqfeng$@$bjtu.edu.cn, klavdija.kutnar@upr.si, dragan.marusic@upr.si, dwyang$@$bjtu.edu.cn}

\author{Yan-Quan Feng$^{\rm a}$, Klavdija Kutnar$^{{\rm b,c}}$,
Dragan Maru\v si\v c$^{{\rm b,c,d}}$, Da-Wei Yang$^{\rm a}$\footnotemark\\
{\small\em $^{\rm a}$Department of Mathematics, Beijing Jiaotong
University, Beijing 100044, China} \\
{\small\em $^{\rm b}$University of Primorska, UP FAMNIT, Glagolja\v ska 8, 6000 Koper, Slovenia} \\
{\small\em $^{\rm c}$University of Primorska, UP IAM, Muzejski trg 2, 6000 Koper, Slovenia}\\
{\small\em $^{\rm d}$University of Ljubljana, UL PEF, Kardeljeva pl. 16, 1000 Ljubljana, Slovenia} \\
}

\date{}
 \maketitle

\begin{abstract}
It was proved in
[Y.-Q. Feng, C. H. Li and J.-X. Zhou, Symmetric cubic graphs with solvable
automorphism groups, {\em European J. Combin.} {\bf 45} (2015), 1-11]
that a cubic symmetric graph with a solvable
automorphism group is either a Cayley graph or a $2$-regular graph of type $2^2$,
that is, a graph with no automorphism of order $2$ interchanging two adjacent vertices.
In this paper an infinite family of non-Cayley cubic $2$-regular graphs of type $2^2$ with a solvable automorphism group is constructed.
 The smallest graph in this family has order 6174.

\bigskip
\f {\bf Keywords:} Symmetric graph, non-Cayley graph, regular cover.\\
{\bf 2010 Mathematics Subject Classification:} 05C25, 20B25.

\end{abstract}

\section{Introduction}
\noindent

Throughout this paper, all groups are finite and all graphs are finite, undirected and simple.
Let $G$ be a permutation group on a set $\O$ and let $\a\in \O$.
We let $G_{\a}$ denote the {\em stabilizer} of
$\a$ in $G$, that is, the subgroup of $G$ fixing the point $\a$.
The group $G$ is {\it semiregular} on
$\O$ if $G_{\a}=1$ for any $\a\in\O$, and {\it regular} if $G$ is transitive and semiregular.
We let $\mz_n$, $\mz_n^*$ and $S_n$ denote
the cyclic group of order $n$,
the multiplicative group of units of $\mz_n$ and the symmetric group of degree $n$, respectively.

For a graph $\G$, we denote its vertex set,
edge set and automorphism group by $V(\G)$, $E(\G)$ and $\Aut(\G)$,
respectively. For a non-negative integer $s$, an {\it $s$-arc} in a graph $\G$ is an ordered $(s+1)$-tuple $(v_0,\ldots,v_s)$ of vertices of $\G$ such that
$v_{i-1}$ is adjacent to $v_i$ for $1\leq i\leq s$, and $v_{i-1}\neq
v_{i+1}$ for $1\leq i<s$. Note that a $0$-arc is just a vertex. A graph $\G$ is {\em $(G,s)$-arc-transitive}
or {\em $(G,s)$-regular} for $G\leq \Aut(\G)$, if $G$ is  transitive or regular on the set of $s$-arcs of $\G$, respectively, and we also say that $G$ is {\em $s$-arc-transitive} or {\em $s$-regular} on $\G$. In particular, $G$ is {\em regular} if it is $0$-regular.
A graph $\G$ is {\em $s$-arc-transitive}
or {\em $s$-regular} if it is $(\Aut(\G),s)$-arc-transitive or $(\Aut(\G),s)$-regular, respectively. Note that $0$-arc-transitive  and $1$-arc-transitive correspond to the terms {\em vertex-transitive} and {\em symmetric}, respectively.

Vertex stabilizers
of connected cubic symmetric graphs were determined in \cite{DM}.
Taking into account the possible isomorphism types for the pair consisting of
a vertex-stabilizer and an edge-stabilizer, the full automorphism groups of
connected cubic symmetric graphs fall into seven classes (see~\cite{CN}).
In particular, for a connected cubic $(G,2)$-regular
graph, if $G$ has an involution flipping an edge, then
it is said to be of {\em type $2^1$}; otherwise, it is of {\em type $2^2$}.
Graphs of type $2^2$ are extremely rare;
there are only nine graphs of type $2^2$ in Conder's list of all cubic
symmetric graphs up to order $10, 000$ \cite{Conder}.

Many cubic symmetric graphs are Cayley graphs, but there are also
examples of non-Cayley graphs among them, such as the
Petersen graph and the Coxeter graph.
For convenience such graphs will be referred to as {\em VNC-graphs}.
Many publications have
investigated VNC-graphs from different perspectives.
For example,
a lot of constructions of VNC-graphs come as a result
of the search for non-Cayley numbers,
that is, numbers for which a VNC-graph of that order exists
(see, for example, \cite{LS,DM83,AMP,MP,MP1,AS}).

The problem of classifying VNC-graphs of small valencies, in particular cubic graphs,
has received a considerable attention (see, for example, \cite{ZF}).
Recently, Feng, Li and Zhou  \cite{FLZ} proved that a connected
cubic symmetric VNC-graph, admitting a solvable arc-transitive group of automorphisms,
is of type $2^2$, and that further such a graph must be a regular cover of the complete bipartite graph $K_{3,3}$  (see Section 3 for the definition of regular covers).
From Conder's list \cite{Conder},
the smallest such graph has order 6174.
In fact, to the best of our knowledge, this graph was the only known
graph of this kind prior to our construction given in this paper.
It is worth mentioning that the family contains a subfamily of
symmetric elementary abelian covers of the Pappus graph of order $18$,
which was overlooked in~\cite{Oh3}.

\section{Preliminaries}
\noindent

Let $G$ be a group and let $M\leq G$.
A subgroup $N$ of $G$ is a {\em normal complement} of $M$
in $G$ if $N\unlhd G$, $N\cap M=1$, and $NM=G$. A normal
complement of a Sylow $p$-subgroup is called the {\em normal $p$-complement} in $G$, that is, a normal $p$-complement in $G$ is a normal Hall $p'$-subgroup of $G$.
The following proposition comes from~\cite[(39.2)]{Aschbacher}.

\begin{prop}\label{p-nilpotent}
Let $G$ be a group. If $p$ is the smallest prime divisor of the order $|G|$
and $G$ has cyclic Sylow $p$-subgroups, then $G$ has a normal $p$-complement.
\end{prop}

From~\cite[Theorem~5.1]{CN}, we have the following proposition.

\begin{prop}\label{prop=2^2}
Let $\G$ be a connected cubic $(G,s)$-regular graph. Then the following hold.

\begin{enumerate}
\itemsep=0pt
\item [\rm (1)] If $G$ has an arc-transitive subgroup of type $2^2$, then
$s=2$ or $3$;
\item [\rm (2)] If $G$ is of type $2^2$, then $G$ has no $1$-regular subgroup.
\end{enumerate}
\end{prop}

Let $\G$ be a graph, and let $K\leq \Aut(\G)$.
The {\em quotient graph} $\G_K$ of $\G$ relative to $K$ is defined as the graph with vertices
the orbits of $K$ on $V(\G)$, with two orbits being
adjacent if there is an edge in
$\G$ between those two orbits.
In view of~\cite[Theorem~9]{L}, we have the following proposition.

\begin{prop}\label{prop=3orbits}
Let $\G$ be a connected cubic $(G,s)$-regular graph for $s\geq 1$, and let $N\unlhd G$. If $N$ has more than two orbits on $V(\G)$, then $N$ is semiregular and the quotient graph $\G_N$ is a cubic $(G/N,s)$-regular graph with $N$ as the kernel
of $G$ acting on $V(\G_N)$.

\end{prop}

 A connected cubic $(G,s)$-regular graph is said to be {\em $G$-basic}, if
$G$ has no non-trivial normal subgroups with more than two orbits on $V(\G)$.
By~\cite[Theorem~1.1]{FLZ}, we have the following proposition.

\begin{prop}\label{prop=basic}
Let $G$ be solvable and let $\G$ be a connected cubic $(G,3)$-regular graph.
If $\G$ is $G$-basic, then $\G\cong K_{3,3}$ and $G\cong S_3^2\rtimes\mz_2$.

\end{prop}

\section{Main result}
\noindent

We first construct connected cubic 2-regular graphs as regular covers of $K_{3,3}$.

\medskip

\begin{figure} [t]
\begin{center}
\unitlength 0.7 cm
\begin{picture} (8,8)

\put (0.5, 7.5){\circle*{0.2}} \put (6.5, 7.5){\circle*{0.2}}
\put (0.5, 4.5){\circle*{0.2}} \put (6.5, 4.5){\circle*{0.2}}
\put (0.5, 1.5){\circle*{0.2}} \put (6.5, 1.5){\circle*{0.2}}

\thicklines \put (0.5,7.5) {\vector(1,0){3}} \put (0.5,7.5) {\line(1,0){6}}
            \put (0.5,4.5) {\vector(1,0){4}} \put (0.5,4.5) {\line(1,0){6}}
            \put (0.5,1.5) {\vector(1,0){3}} \put (0.5,1.5) {\line(1,0){6}}
            \put (0.5,4.5) {\vector(2,1){4}} \put (0.5,4.5) {\line(2,1){6}}
            \put (0.5,4.5) {\vector(2,-1){4}} \put (0.5,4.5) {\line(2,-1){6}}

\thinlines  \put (0.5,7.5) {\vector(2,-1){1.5}} \put (0.5,7.5) {\line(2,-1){6}}
            \put (0.5,7.5) {\vector(1,-1){1.4}} \put (0.5,7.5) {\line(1,-1){6}}
            \put (0.5,1.5) {\vector(1,1){1.4}} \put (0.5,1.5) {\line(1,1){6}}
            \put (0.5,1.5) {\vector(2,1){1.5}} \put (0.5,1.5) {\line(2,1){6}}

\put (-0.2, 7.4){{\bf v}} \put (-0.2, 4.4){{\bf u}} \put (-0.2, 1.4){{\bf w}}
\put (6.8, 7.4){{\bf y}} \put (6.8, 4.4){{\bf x}} \put (6.8, 1.4){{\bf z}}

\put (3.1, 7.7){$1$}   \put (3.1, 0.85){$1$} \put (4.2, 6.6){$1$}
\put (4.2, 4.7){$1$}   \put (4.2, 2){$1$}    \put (1.9, 6.9){$h^{-1}b$}
\put (1.2, 5.8){$h$} \put (0.9, 2.6){$hc$}
\put (2, 1.7){$h^{-1}a$}
\end{picture}
\end{center}\vspace{-1cm}
\caption {\label{fig:K33}The complete bipartite graph $K_{3,3}$ with a voltage assignment $\phi$.}\label{fig=1}
\end{figure}
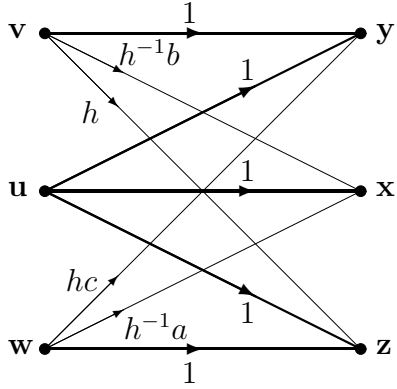

\f {\bf Construction:} Let $n$ be an integer such that $n\geq7$ and
that the equation $x^2+x+1=0$ has a solution $r$ in $\mz_n$. Then $r$ is an element of order 3 in $\mz_n^*$, and,
by~\cite[Lemma 3.3]{FL}, the prime decomposition of $n$ is $3^tq_1^{e_1}\cdots q_s^{e_s}$ with
$t\leq1$, $s\geq 1$, $e_i\geq1$ and $3\di (q_i-1)$ for $1\leq i\leq s$. In particular, $n$ is odd.
Let $K=\lg a,b,c,h~|~a^n=b^n=c^n=h^3=[a,b]=[a,c]=[b,c]=1,a^h=a^r,b^h=b^r,c^h=c^r\rg$
be a group. Then $K\cong\mz_n^3\rtimes\mz_3$ and $K$ has odd order $3n^3$.
Denote by $V(K_{3,3})=\{{\bf u}, {\bf v}, {\bf w}, {\bf x}, {\bf y}, {\bf z}\}$
the vertex set of $K_{3,3}$ such that the vertices from the set $\{{\bf u}, {\bf v}, {\bf w}\}$
are adjacent to the vertices from the set $\{{\bf x}, {\bf y}, {\bf z}\}$, see Figure~\ref{fig:K33}.
The graph $\mathcal{NCG}_{18n^3}$ is defined to have the vertex set
$V(\mathcal{NCG}_{18n^3})=V(K_{3,3})\times K$ and edge set
$$\begin{array}{lll}
& E(\mathcal{NCG}_{18n^3})=\{\{({\bf u},g), ({\bf x},g)\}, \{({\bf u},g),({\bf y},g)\},
\{({\bf u},g),({\bf z},g)\}, \{({\bf v}, g),({\bf x}, gh^{-1}b)\},\\
&~~~~~~~~~~~~~~~~~~~~~~ \{({\bf v},g),({\bf y},g)\}, \{({\bf v},g),({\bf z},gh)\}, \{({\bf w},g),({\bf x},gh^{-1}a)\},\{({\bf w},g),({\bf y},ghc)\},\\
&~~~~~~~~~~~~~~~~~~~~~~ \{({\bf w},g),({\bf z},g)\}~|~ g\in K\}.
\end{array}$$
\f Clearly, $\mathcal{NCG}_{18n^3}$ is a bipartite graph.

\begin{theorem}\label{the=main}
The graph $\mathcal{NCG}_{18n^3}$ is a connected cubic symmetric non-Cayley
graph and its automorphism group is solvable and of type $2^2$.
\end{theorem}

\medskip

To prove Theorem~\ref{the=main}, we need to introduce the voltage graph.
Let $\G$ be a connected graph and $K$ a group. Assign to each arc
$(u, v)$ of $\G$ a voltage $\phi(u, v)$ such that $\phi(u, v)\in K$ and
$\phi(u, v)=\phi(v, u)^{-1}$, where $\phi: \G\mapsto K$
is a {\em voltage assignment} of $\G$. Let $\G\times_{\phi} K$
be the {\em voltage graph} obtained from $\phi$ in the following way:
$V(\G)\times K$ is its vertex set and $\{\{(u,a),(v,a\phi(u,v))\}~|~\{u,v\}\in E(\G),a\in K\}$ is its edge set.
Now, $\G\times_{\phi} K$ is a {\em regular cover} (or a {\em $K$-cover}) of $\G$, and
the graph $\G$ is a {\em base graph}.
Moreover, the quotient graph $(\G\times_{\phi} K)_K$ is isomorphic to $\G$,
and $\G\times_{\phi}K$ is connected if and only if the voltages on the arcs generate the voltage group $K$. The projection onto the first
coordinate $\pi: \G\times_{\phi} K\mapsto \G$ is a {\em regular
$K$-covering projection}, where the group $K$ acts semiregularly via a
left multiplication on itself. Clearly, $\pi$ induces an isomorphism $\bar{\pi}$ from the quotient graph $(\G\times_{\phi} K)_K$ to $\G$.
We say that an automorphism $\a$ of $\G$ {\em lifts} to an automorphism $\widetilde{\a}$ of
$\G\times_{\phi} K$ if $\widetilde{\a}\pi=\pi\a$. In this case, $\widetilde{\a}$ is a {\em lift} of $\a$. In particular, $K$ is the lift of the identify group of $\Aut(\G)$, and if an automorphism $\a\in \Aut(\G)$ has a lift $\widetilde{\a}$,
then $K\widetilde{\a}$ are all lifts of $\a$.
Let $F=N_{{\rm Aut}(\G\times_{\phi} K)}(K)$, the normalizer of $K$ in $\Aut(\G\times_{\phi} K)$, and let $L$ be the largest subgroup of $\Aut(\G)$ that can be lifted.  Since $K\unlhd F$, each automorphism $\widetilde{\a}$ in $F$ induces an automorphism of $(\G\times_{\phi} K)_K$ and hence an automorphism $\a$ of $\G$ via $\bar{\pi}$, so
$\widetilde{\a}$ is a lift of $\a$.
On the other hand, it is easy to check that the lifts of each automorphism in $L$
map orbits  to orbits of $K$, and so normalize $K$. Therefore, there exists an epimorphism $\psi$
from $F$ to $L$ with kernel $K$, and thus $F/K\cong L$.
For an extensive treatment of regular covering we refer the reader to \cite{23}.

The problem whether an automorphism $\a$ of $\G$ lifts
can be grasped in terms of voltages as follows.
Observe that a voltage assignment on arcs extends to a voltage
assignment on walks in a natural way. For $\alpha\in \Aut(\G)$,
we define a function $\bar{\alpha}$ from the set of voltages on
fundamental closed walks based at a fixed vertex $v$ in $V(\G)$ to
the voltage group $K$ by $(\phi(C))^{\bar{\alpha}}=\phi(C^{\alpha})$,
where $C$ ranges over all fundamental closed walks at $v$, and
$\phi(C)$ and $\phi(C^{\alpha})$ are the voltages of $C$ and $C^{\alpha}$,
respectively. The next
proposition is a special case of \cite[Theorem 4.2]{23}.

\begin{prop}\label{lifts}
If ${\pi}: \G\times_{\phi} K \mapsto \G$ is a connected regular $K$-covering,
then an automorphism $\alpha$ of $\G$ lifts if and only if $\bar{\alpha}$ extends to an automorphism of $K$.
\end{prop}

Now, we are ready to prove Theorem~\ref{the=main}.

\medskip

\f {\bf Proof of Theorem~\ref{the=main}:} By Construction, $\mathcal{NCG}_{18n^3}$
is the voltage graph $K_{3,3}\times_\phi K$ with the voltage assignment $\phi$ as depicted in Figure~\ref{fig=1}. Let $\mathcal{P}: \mathcal{NCG}_{18n^3}\mapsto K_{3,3}$
be the covering projection.
Let $A=\Aut(\mathcal{NCG}_{18n^3})$ and $F=N_{A}(K)$.
Then $L=F/K$ is the largest subgroup of $\Aut(K_{3,3})$,
which can be lifted along $\mathcal{P}$.

Denote by $i_1i_2\cdots i_s$ the cycle having the consecutively adjacent vertices $i_1,
i_2,\ldots,i_s$. There are four fundamental closed walks based at the vertex {\bf u}
in $K_{3,3}$, that is, {\bf uyvz}, {\bf uzwx}, {\bf uyvx} and {\bf uzwy},
which are generated by the four cotree arcs ({\bf v}, {\bf z}), ({\bf w}, {\bf x}),
({\bf v}, {\bf x}), and ({\bf w}, {\bf y}), respectively.
Since $\lg \phi({\bf v}, {\bf z}), \phi({\bf w}, {\bf x}),
\phi({\bf v}, {\bf x}),\phi({\bf w}, {\bf y})\rg=\lg h,h^{-1}a,h^{-1}b,hc\rg=K$,
the graph $\mathcal{NCG}_{18n^3}$ is connected.

Define four permutations on $V(K_{3,3})$ as follows:
$$\a_1=({\bf uvw}),~\a_2=({\bf xyz}),~\b=({\bf ux})({\bf vy})({\bf wz}),~\d=({\bf vywz})({\bf ux}).$$
It is easy to check that $\Aut(K_{3,3})=\lg \a_1,\a_2,\b,\d\rg\cong (S_3\times S_3)\rtimes\mz_2$, and $\lg \a_1,\a_2,\d\rg\cong(\mz_3\times\mz_3)\rtimes\mz_4$ is 2-regular. Clearly, each involution in $\lg \a_1,\a_2,\d\rg$ fixes the bipartite parts of $K_{3,3}$, and then the subgroup $\lg \a_1,\a_2,\d\rg$ is of type $2^2$ and contains no regular subgroup. Since $\a_1^\d=\a_2$ and $\a_2^\d=\a_1^{-1}$, we have that $\lg \a_1,\a_2,\d\rg$ has no normal subgroup of order $3$.

Under $\a_1$,
$\a_2$, $\b$ and $\d$, each walk of $K_{3,3}$
is mapped to a walk of the same length. We list all these walks and their voltages in Table~\ref{table=1},
in which $C$ denotes a fundamental closed
walk of $K_{3,3}$ based at the vertex ${\bf u}$ and
$\phi(C)$ denotes the voltage on $C$.

\begin{table}[ht]

\begin{center}
\begin{tabular}{|l|l|l|l|l|l|}

\hline
$C$ &$\phi(C)$ &$C^{\a_1}$&$\phi(C^{\a_1})$       &$C^{\a_2}$&$\phi(C^{\a_2})$ \\
  \hline
{\bf uyvz} & $h$        &{\bf vywz}  &$hc^{-r}$                &{\bf uzvx}  &$hb$ \\
{\bf uzwx} & $h^{-1}a$  &{\bf vzux}  &$h^{-1}b^{-r}$           &{\bf uxwy}  &$h^{-1}a^{-r^2}c$\\
{\bf uyvx} & $h^{-1}b$  &{\bf vywx}  &$h^{-1}a^rb^{-r}c^{-r^2}$&{\bf uzvy}  &$h^{-1}$\\
{\bf uzwy} & $hc$       &{\bf vzuy}  &$h$                      &{\bf uxwz}  &$ha^{-r}$  \\
  \hline
\hline
$C$ &$\phi(C)$ &$C^{\b}$ &$\phi(C^{\b})$         &$C^{\d}$ &$\phi(C^{\d})$ \\
  \hline
{\bf uyvz} & $h$        &{\bf xvyw} &$h^{-1}ab^{-r^2}c^{-r}$ &{\bf xwyv} &$ha^{-r}bc^{r^2}$ \\
{\bf uzwx}& $h^{-1}a$  &{\bf xwzu} &$ha^{-r}$                &{\bf xvzu} &$h^{-1}b^{-r^2}$\\
{\bf uyvx} & $h^{-1}b$  &{\bf xvyu} &$hb^{-r}$                &{\bf xwyu} &$h^{-1}a^{-r^2}c$\\
{\bf uzwy} & $hc$       &{\bf xwzv} &$h^{-1}a^{-r^2}b$        &{\bf xvzw} &$hab^{-r}$  \\
 \hline
\end{tabular}
\end{center}
\vskip -0.5cm
\caption{Fundamental walks and their images with corresponding voltages.} \label{table=1}
\end{table}

Let $\bar{\a_1}$ be the map defined by $\phi(C)^{\bar{\a_1}}=\phi(C^{\a_1})$,
where $C$ ranges over the four fundamental
closed walks of $K_{3,3}$ based at the vertex {\bf u}. Similarly, we can define $\bar{\a_2}$, $\bar{\b}$ and $\bar{\d}$.
Recall that $r^2+r+1=0~(\mod n)$. The following equations will be used frequently:
$$[a,b]=[a,c]=[b,c]=1,\ ah=ha^r,\ bh=hb^r,\ ch=hc^r.$$
By Table~\ref{table=1}, one may easily check that $\bar{\a_1}$, $\bar{\a_2}$ and $\bar{\d}$
extend to three automorphisms of $K$ induced by
$a\mapsto b^{-r}c^{-1}$, $b\mapsto a^rb^{-r}c^r$, $c\mapsto c^r$, $h\mapsto hc^{-r}$;
$a\mapsto a^{-r^2}b^{r^2}c$, $b\mapsto b^{r^2}$, $c\mapsto a^{-r}b^{-1}$, $h\mapsto hb$;
and $a\mapsto a^{-1}c^{r}$, $b\mapsto a^rb^{r^2}c^{-r^2}$, $c\mapsto a^{-r^2}b^{r^2}c^{-r^2}$, $h\mapsto ha^{-r}bc^{r^2}$, respectively.
Hence, by Proposition~\ref{lifts}, $\a_{1}$, $\a_2$ and $\d$ lift.

Suppose $\bar{\b}$ extends to an automorphism of $K$, say $\b^*$.
By Table~\ref{table=1}, $h^{\b^*}=h^{-1}ab^{-r^2}c^{-r}$
and $(h^{-1}a)^{\b^*}=ha^{-r}$. Thus $a^{\b^*}=(h\cdot h^{-1}a)^{\b^*}=
h^{\b^*}\cdot (h^{-1}a)^{\b^*}=h^{-1}ab^{-r^2}c^{-r}\cdot ha^{-r}=b^{-1}c^{-r^2}$
and $(a^{\b^*})^{h^{\b^*}}=(b^{-1}c^{-r^2})^{h^{-1}ab^{-r^2}c^{-r}}=b^{-r^2}c^{-r}$.
Since $a^h=a^r$, we have $(a^{\b^*})^{h^{\b^*}}=(a^h)^{\b^*}=(a^{\b^*})^r$, that is,
$b^{-r^2}c^{-r}=(b^{-1}c^{-r^2})^r=b^{-r}c^{-1}$.
Hence $-r=-1$ and $r=1$. It follows that $3=0~(\mod n)$ because $r^2+r+1=0~(\mod n)$,
contradicting $n\geq7$.

By Proposition~\ref{lifts}, we can conclude that $\b$ does not lift. Since $\a_1,\a_2,
\d$ lift and $|\Aut(K_{3,3}):\lg\a_1, \a_2,\d\rg |=2$, the largest lifted group is $L=\lg \a_1,\a_2,\d\rg$.
Since $F/K=L$ and $L$ is 2-regular, $\mathcal{NCG}_{18n^3}$ is $(F,2)$-regular and $F$ is solvable.

Suppose $F$ is of type $2^1$. Then $F$ has an involution $g$ reversing an edge in $\mathcal{NCG}_{18n^3}$, and thus $gK$ is an involution in $F/K$ reversing an edge in $(\mathcal{NCG}_{18n^3})_K=K_{3,3}$,
which is impossible because $L$ is of type $2^2$. Thus $F$ is of type $2^2$.

To complete the proof we need to show that $F=A$ and that
$\mathcal{NCG}_{18n^3}$ is not a Cayley graph. We first prove that $F$ has no regular subgroup.
Suppose, on the contrary, that $F$ has a regular subgroup $R$ on $V(\mathcal{NCG}_{18n^3})$.
Then $R$ has order twice an odd integer. Since $\mathcal{NCG}_{18n^3}$ is  bipartite, $R$ has an involution $g$ interchanging the two bipartite parts of $\mathcal{NCG}_{18n^3}$.
Since $F/K\cong \lg\a_1,\a_2,\d\rg\cong(\mz_3\times\mz_3)\rtimes\mz_4$
and $K$ has the odd order $3n^3$, the Sylow 2-subgroups of $F$ are isomorphic to $\mz_4$.
By Proposition~\ref{p-nilpotent}, $F$ has a normal Hall $2'$-subgroup $H$, which has two orbits as the two bipartite parts of $\mathcal{NCG}_{18n^3}$
with vertex stabilizer isomorphic to $\mz_3$.
Then $H\lg g\rg$ is a 1-regular subgroup of $F$, contrary to
Proposition~\ref{prop=2^2}~(2). Hence $F$ has no regular subgroup on $V(\mathcal{NCG}_{18n^3})$.

Finally, suppose that $A\neq F$. Since $A$ has an arc-transitive proper subgroup $F$
of type $2^2$, by Proposition~\ref{prop=2^2}~(1), the group $A$ is $3$-regular.
This implies that $|A:F|=2$ and $F\unlhd A$. Since $A/F\cong\mz_2$ and $F$ is solvable, $A$ is solvable.

Let $H$ be a maximal normal subgroup of $A$ having
at least three orbits on $V(\mathcal{NCG}_{18n^3})$. By Proposition~\ref{prop=3orbits},
$H$ is semiregular and the quotient graph $(\mathcal{NCG}_{18n^3})_H$ is
$(A/H,3)$-regular. By the maximality of $H$, $(\mathcal{NCG}_{18n^3})_H$ is
$A/H$-basic. Since $A/H$ is solvable, Proposition~\ref{prop=basic} implies that  $(\mathcal{NCG}_{18n^3})_H\cong K_{3,3}$
and $A/H\cong S_3^2\rtimes \mz_2$. It follows that $|H|=|V(\mathcal{NCG}_{18n^3})|/6=|K|$.
Since $|A:F|=2$ and $|H/H\cap F|=|HF/F|\di |A/F|$, we have $|H/H\cap F|=1$ or $2$, and since $|H|$ is odd, $|H/H\cap F|=1$, that is , $H\leq F$.

Recall that $F=N_A(K)$. Since $A\neq F$, the subgroup $K$ is not normal in $A$, and thus $H\neq K$.
Since $H\unlhd F$, we have $H\cap K\unlhd F$ and $1\neq HK/K\unlhd F/K$.
Since $|K|=|H|$ is odd and $L=F/K\cong \mz_3^2\rtimes\mz_4$, the quotient group $HK/K$ is a non-trivial 3-group, and since $L$ has no normal subgroup of order $3$, we have $HK/K\cong\mz_3^2$. It follows that
$|H\cap K|=\frac{|H||K|}{|HK|}=\frac{1}{9}|K|$
and thus $|H\cap K|=\frac{1}{54}|V(\mathcal{NCG}_{18n^3})|$ since $|V(\mathcal{NCG}_{18n^3})|=6|K|$.
Since $F$ is 2-regular on $\mathcal{NCG}_{18n^3}$ and $H\cap K\unlhd F$, Proposition~\ref{prop=3orbits}
implies that
the quotient graph $(\mathcal{NCG}_{18n^3})_{H\cap K}$
is a connected cubic $(F/H\cap K,2)$-regular graph of order
$54$. Moreover, since $F$ has no regular subgroup,
$F/H\cap K$ has no regular subgroup on $V((\mathcal{NCG}_{18n^3})_{H\cap K})$.
However, by~\cite{Conder} there is only one connected cubic symmetric graph
of order 54, the graph $F_{54}$, which is 2-regular, of girth 6, and, by~\cite[Theorem~1.1]{KM},
a Cayley graph. It follows that $(\mathcal{NCG}_{18n^3})_{H\cap K}\cong F_{54}$ and $F/H\cap K=\Aut((\mathcal{NCG}_{18n^3})_{H\cap K})$,
so that $F/H\cap K$ has a regular subgroup, a contradiction.

Thus $A=F$  is solvable, of type $2^2$ and has no regular subgroup.
In other words, $\G$ is non-Cayley as claimed. This completes the proof. \hfill\qed

\f {\bf Remark:} By the proof of Theorem~\ref{the=main}, $K\unlhd \Aut(\mathcal{NCG}_{18n^3})$.
When $n$ is a prime $p$, that is, $K\cong\mz_p^3\rtimes\mz_3$
with $3\di (p-1)$, the group $K$ has a characteristic Sylow $p$-subgroup $P$
and $P\cong\mz_p^3$. Thus $P\unlhd \Aut(\mathcal{NCG}_{18p^3})$.
Clearly, $P$ has more than two orbits
on $V(\mathcal{NCG}_{18p^3})$. By Proposition~\ref{prop=3orbits}, $(\mathcal{NCG}_{18p^3})_P$ is a connected cubic symmetric graph of order $18$. By~\cite{Conder},
up to isomorphism, there is only one connected
cubic symmetric graph of order 18, that is, the Pappus graph.
Hence $\mathcal{NCG}_{18p^3}$ is a connected 2-regular $\mz_p^3$-cover of the Pappus graph
with $3\di (p-1)$. These graphs were overlooked in~\cite[Theorem~3.1]{Oh3}.

\medskip
\f {\bf Acknowledgement:} This work was supported by the National Natural Science Foundation of China (11571035, 11231008) and by the 111 Project of China (B16002).
The work of K.~K. was supported in part by the Slovenian Research Agency
(research program P1-0285 and research projects N1-0032, N1-0038, J1-6720, J1-6743,
and J1-7051), in part by WoodWisdom-Net+, W$^3$B, and in part by NSFC project 11561021.
The work of D.~M. was supported in part by the Slovenian Research Agency (I0-0035, research program P1-0285 and research projects N1-0032, N1-0038, J1-5433,  J1-6720,
and J1-7051), and in part by H2020 Teaming InnoRenew CoE.

\end{document}